\documentclass[12pt,reqno,a4paper]{amsart}
\usepackage[latin1]{inputenc}

\usepackage{amsmath}
\usepackage{amsfonts}
\usepackage{amssymb}
\usepackage{listings}
\usepackage[colorlinks=true,linkcolor=blue,citecolor=blue,urlcolor=blue]{hyperref}

\title{Multivariate Polynomials in Sage}
\author{Viviane Pons}
\address{Institut Gaspard Monge, Universit\'e Paris-Est,
Marne-la-Vall\'ee, France}
\date{}

\begin{document}

\begin{abstract}
We have developed a patch implementing multivariate polynomials seen as a multi-base algebra. The patch is to be released into the software Sage and can already be found within the Sage-Combinat distribution. One can use our patch to define a polynomial in a set of indexed variables and expand it into a linear basis of the multivariate polynomials. So far, we have the Schubert polynomials, the Key polynomials of types $A$, $B$, $C$, or $D$, the Grothendieck polynomials and the non-symmetric Macdonald polynomials. One can also use a double set of variables and work with specific double-linear bases like the double Schubert polynomials or double Grothendieck polynomials. Our implementation is based on a definition of the basis using \textit{divided difference} operators and one can also define new bases using these operators. 
\end{abstract}

\maketitle

\thispagestyle{myheadings}
\font\rms=cmr8
\font\its=cmti8
\font\bfs=cmbx8

\markright{\its S\'eminaire Lotharingien de
Combinatoire \bfs 66 \rms (2011), Article~B66z\hfill}
\def\thepage{}

\section*{Introduction}

Multivariate polynomials and their bases appear in many combinatorial problems and one often needs to define a polynomial as a formal sum of elements that live in a specified basis. The usual implementation of multivariate polynomials is done as a tensor product of polynomials in one variable. But one can not consider the variables all together: all bases have to be defined by a product of bases of polynomials in one variable. It appeared to us that a clear and handy implementation of multivariate polynomials from a combinatorial point of view would be useful not only to our work but to the community. 

Our approach is based on \textit{divided difference} operators and their interpretation in terms of linear bases of the multivariate polynomials algebra as explained by Lascoux \cite{DUMMIES}. We define simple operators of types $A$, $B$, $C$, and $D$ in Section \ref{OPER} and use them to create the \textit{divided difference} operators in Section \ref{DIFFDIV}. We then explain in Section \ref{BASES} how these operators allow us to define linear bases of the multivariate polynomials.

Our software has been implemented in Sage, we explain our choice in Section \ref{SAGE}. The full description of the implemented features with code examples can be found in Section \ref{SOFTWARE}. The development process is not finished at the time of submission of this paper and what the program now needs most is to be tested by many users so that bugs can be reported and suggestions be made on how to improve the features.

\section{Choosing Sage}
\label{SAGE}

Sage is a mathematics software created in 2005. It is free and open source which is the main reason why we have chosen it. It is developed by many researchers around the world and one can join the very lively community who works specifically on combinatorics within the Sage-Combinat project. The Sage-Combinat community has been created in 2001 and was previously known as Mupad-Combinat. It has moved to Sage in 2008 keeping its main purpose: developing specific tools for algebraic combinatorics and sharing the programs among researchers.

From the beginning, we wanted not only to develop in Sage but to be part of the main project by adding our implementation to the software. The program is still in test mode within the Sage-combinat distribution (see the installation process in Section \ref{INSTALL}).

Working in Sage also allowed us to use previous work and structures developed by the community. We worked a lot with Weyl Groups which had already been implemented in Sage. We also used the standard implementation of multi-base algebras as a base for our own work. 

And finally, we hope that being part of Sage will help our software to spread within the community and thus become useful to the largest possible community. 

\pagenumbering{arabic}
\addtocounter{page}{1}
\markboth{\SMALL VIVIANE PONS}{\SMALL MULTIVARIABLE POLYNOMIALS IN SAGE}

\section{Multi-base polynomials}
\subsection{Type $A$, $B$, $C$, $D$ operators}
\label{OPER}

At first, we need to define simple operations on vectors of integers. Let $v \in \mathbb{Z}^n$, we have the following operators corresponding to the root system of respective types $A$, $B$, $C$, $D$:
\begin{align}
vs_i &= (\ldots, v_{i+1},v_i, \ldots) &\text{ for }1 \leq i < n \\
vs_i^B = vs_i^C &= (\ldots, -v_i, \ldots) &\text{ for } 1 \leq i \leq n \\
vs_i^D &= (\ldots, -v_i, -v_{i-1}, \ldots) &\text{ for } 2 \leq i \leq n
\end{align}
The group generated by  $s_1, \ldots, s_{n-1}$ (respectively  $s_1, \ldots, s_{n-1},s_n^B$, and $s_1, \ldots,\\ s_{n-1},s_n^D$) is the Weyl group of type $A$ (respectively $B$ or $C$, $D$). The operators satisfy the \textit{braid relations}:
\begin{align}
s_is_{i+1}s_i &= s_{i+1}s_is_{i+1} &\text{ and } s_is_j = s_js_i \text{, } |i-j| \neq 1 \\
s_{n-1}s_n^Bs_{n-1}s_n^B &= s_n^Bs_{n-1}s_n^Bs_{n-1} &\text{ and } s_is_n^B = s_n^Bs_i \text{, } i \leq n-2 \\
s_{n-2}s_n^Ds_{n-2} &= s_n^Ds_{n-2}s_n^D &\text{ and } s_is_n^D = s_n^Ds_i \text{, } i \neq n-2
\end{align}
The orbit of the vector $[1, 2, \dots , n]$ consists of all permutations of $1, 2,\dots , n$, i.e., $S_n$, for type $A$, all signed permutations for type $B$, $C$, and all signed permutations with an even number of minus signs for type $D$. The elements of the different groups can be denoted by these objects. In the same way, elements of these groups can be seen as a product of operators $s_i$, called a \textit{decomposition}. When the product is of minimal length, it is called a \textit{reduced decomposition}. 

\subsection{Action on polynomials}
\label{DIFFDIV}

We now have a natural action of the Weyl groups on polynomials. Indeed, let $x = (x_1, x_2, \ldots, x_n)$  be a set of variables and for $v \in \mathbb{Z}^n$, let $x^v$ stand for the monomial
\begin{equation}
x_1^{v_1}x_2^{v_2}\ldots x_n^{v_n}.
\end{equation}
A polynomial in the variables $x_1, \ldots, x_n$ can therefore be seen as a formal sum of vectors and the action of the operator $s_i$ becomes an action on polynomials:
\begin{equation}
x^vs_i = x^{vs_i}.
\end{equation}
From these simple operators $s_i$, we can now define the \textit{divided difference} operators. For type $A$, we have:
\begin{align}
f\partial_i &:= \frac{f - f^{s_i}}{x_i - x_{i+1}} \\
f\pi_i &:= \frac{x_if - x_{i+1}f^{s_i}}{x_i - x_{i+1}} = f.(x_i\partial_i) \\
f\hat{\pi}_i &:=\frac{(f-f^{s_i})x_{i+1}}{x_i - x_{i+1}} = f.(\pi_i - 1) \\
fT_i &:= f.(\pi_i(t_1 + t_2) - s_it_2)
\end{align}
for $ 1 \leq i \leq n-1$. $\partial_i$ is the \textit{Newton divided difference}, $\pi_i$ and $\hat{\pi}_i$ are the \textit{isobaric divided differences} and $T_i$ is the generator of the \textit{Hecke algebra} $\mathcal{H}_2$. As the $s_i$ satisfy the braid relations, all the above operators do. 

If $v_i > v_{i+1}$, the \textit{Newton divided difference} $\partial_i$ can be seen as an operator decrementing the vector degree and summing over all the intermediate monomials between $x^{(\ldots, v_i-1, v_{i+1}, \ldots)}$ and $x^{(\ldots, v_{i+1}, v_{i}-1, \ldots)}$. For example,
\begin{equation}
x^{(4,1)}\partial_1 = x^{(3,1)} + x^{(2,2)} + x^{(1,3)}.
\end{equation}
When $v_i < v_{i+1}$, one just needs to switch $v_i$ and $v_{i+1}$, multiply by $-1$ and do the previous operation. When $v_i =v_{i+1}$, the result is $0$. This description can be used to give a more general, type-free definition. We can see the vectors indexing the monomials as elements of the ambient space of the root system of type $A_{n-1}$. The above $\partial_i$ operation can then be seen as a formal sum of vectors, adding factors of the simple root $(\ldots, -1, 1 \ldots)$ to the original vector. The sign of $v_i - v_{i+1}$ is given by the scalar product between the vector and the $i^{th}$ simple coroot of the ambient space. This definition is equivalent to the previous one. We can use it to define our divided differences in types $B$, $C$, and $D$, using the root systems of respective types $B_n$, $C_n$, and $D_n$. Compared to type $A$, we add the $n^{th}$ simple root and coroot whose definition depends on the type and create the $n^{th}$ divided difference operator:
\begin{align}
\partial_n^B &= \frac{1 - s_n^B}{x_n^{\frac{1}{2}} - x_{n}^{-\frac{1}{2}}} \\
\partial_n^C &= \frac{1 - s_n^C}{x_n - x_n^{-1}} \\
\partial_n^D &= \frac{1 - s_n^D}{x_{n-1}^{-1} - x_n}
\end{align}
The same construction can be done with the isobaric divided differences $\pi$ and $\hat{\pi}$. 

\subsection{Linear bases}
\label{BASES}

We can now use these operators to define linear bases of the ring of multivariate polynomials. Let $(x_1, x_2, \ldots, x_n)$ and $(y_1, y_2, \ldots, y_n)$ be two sets of variables and $\lambda$ a partition of length $n$, i.e., $\lambda_1 \geq \lambda_2 \geq \ldots \geq \lambda_n$. We then define \textit{dominant Schubert polynomials} (respectively \textit{Grothendieck polynomials} and \textit{Key polynomials}) by
\begin{align}
Y_\lambda &:= \prod_{i = 1}^n \prod_{j = 1}^{\lambda_i} (x_i - y_j), \\
G_\lambda &:= \prod_{i = 1}^n \prod_{j = 1}^{\lambda_i} (1-y_jx_i^{-1}), \\
K_\lambda = \hat{K}_\lambda &:= x^\lambda.
\end{align}
We define Schubert polynomials to be all the non-zero images of the dominant Schubert polynomials under products of $\partial_i$ and Grothendieck polynomials to be all the images of the dominant Grothendieck polynomials under products of $\pi_i$. Similarly, the two types of Key polynomials are defined by taking all the images under products of $\pi_i$ or of $\hat{\pi}_i$ respectively. Since the operators satisfy relations, we cannot index the polynomials by the choice of the starting point and the sequence of operators used. Rather, we use weight vectors $v \in \mathbb{N}^n$, the recursive definition being
\begin{align}
Y_{\ldots, v_{i+1}, v_{i}-1, \ldots} = Y_v \partial_i\\
G_{\ldots, v_{i+1}, v_{i}-1, \ldots} = G_v \pi_i\\
K_v\pi_i = K_{vs_i}\\
\hat{K}_v \hat{\pi}_i = \hat{K}_{vs_i},
\end{align}
the inal vectors $v$ satisfying $v_i > v_{i+1}$.

As the operators satisfy braids relations, the order one chooses to apply the recursive rule on a vector does not change the result. There are dominant polynomials in the images of a dominant polynomial in the
Schubert and Grothendieck case; therefore, one has to check consistency, but this is easy. These families constitute triangular bases of the polynomials in $(x_1, \ldots, x_n)$. And one can easily express an arbitrary polynomial in these bases by inverting a triangular matrix. When working with a single set of variables, one can specialize the $y_i$'s to 0 (respectively to 1) and obtain simple Schubert polynomials (respectively Grothendieck polynomials), so that dominant polynomials become
\begin{align}
Y_\lambda &= x^\lambda, \\
G_\lambda &= \prod_{i = 1}^n (1 - x_i^{-1})^{\lambda_i}.
\end{align}
To work with positive exponents on the Grothendieck basis, one can also set $x_i = 1 - x_i^{-1}$. Both versions of the bases are available in our implementation. 

Using the same method, we can also define non-symmetric Macdonald polynomials. In this case, there will be only one generator polynomial, i.e., $M_{0,0 \ldots, 0} = 1$, and we will use both the $T_i$ operator and a raising operator to increase the polynomial degree. The recursive rule is due to Sahi and Knop \cite{SahiKnop} and one can find its description in Lascoux \cite{DUMMIES}.

One can also define type $B$, $C$, and $D$ Key polynomials using the operators defined from the Weyl group of the given types as explained in Section \ref{DIFFDIV}. The corresponding families of polynomials will then be indexed by vectors in $\mathbb{Z}^n$ instead of $\mathbb{N}^n$ and become bases of the Laurent polynomials of $(x_1, \ldots, x_n)$, i.e., with both positive and negative exponents.

\section{Software description}
\label{SOFTWARE}
\subsection{Installation process}
\label{INSTALL}

Our software has been developed as a part of the Sage project. Nevertheless, as it is still in test mode, at the publication time of this paper, it is not yet available on the main Sage distribution, but one can use it within the Sage-Combinat distribution.

\medskip
\noindent
\textsc{Step 1}
\medskip

If Sage is not already installed on a computer, please follow the instructions on the Sage website \cite{SAGE_WEBSITE} to get the installation process corresponding to your system.

To get the latest version of our program, make sure you have the last version of Sage installed. You can upgrade your Sage version by running the following command inside your sage directory:

{\tt
\begin{lstlisting}
./sage -upgrade
\end{lstlisting}
}

The examples below work with Sage 4.7 and later versions.

\medskip
\noindent
\textsc{Step 2}
\medskip

Install Sage-Combinat \cite{SAGE_COMBINAT} by running the following command inside the sage directory:

{\tt \begin{lstlisting}
./sage -combinat install
\end{lstlisting}}

If you encounter problems, you find more information by visiting the Sage-Combinat website \cite{SAGE_COMBINAT}.

\subsection{Define a polynomial}

Sage programming is object-oriented. The object containing the main software methods is called \textit{AbstractPolynomialRing}. One first needs to create this object:

{\tt \begin{lstlisting}
sage: A = AbstractPolynomialRing(QQ)
sage: A
The abstract ring of multivariate polynomials on x over Rational Field
\end{lstlisting}}

\noindent
Here,
\textit{A} represents the abstract algebra. To create an actual polynomial, we need a concrete basis.

{\tt \begin{lstlisting}
sage: m = A.monomial_basis(); m
The ring of multivariate polynomials on x over Rational Field on the monomial basis
\end{lstlisting}}

\noindent
\textit{m} is a concrete basis and we shall use it to create polynomials. Both of the syntaxes presented below can be used. 

{\tt \begin{lstlisting}
sage: pol = m[1,1,2] + m[2,3]; pol
x[1, 1, 2] + x[2, 3, 0]
sage: pol = m([1,1,2]) + m([2,3]); pol
x[1, 1, 2] + x[2, 3, 0]
\end{lstlisting}}

\noindent
$x[1,1,2]$ means $x^{(1,1,2)}=x_1^1x_2^1x_3^2$. One does not have to declare beforehand how many variables are to be used, it will be computed from the size of the vectors. To know on how many variables a polynomial is defined, one can look at its parent. It can also be changed if needed.

{\tt \begin{lstlisting}
sage: pol.parent()
The ring of multivariate polynomials on x over Rational Field with 3 variables on the monomial basis
sage: pol = pol.change_nb_variables(4)
sage: pol
x[1, 1, 2, 0] + x[2, 3, 0, 0]
sage: pol.parent()
The ring of multivariate polynomials on x over Rational Field with 4 variables on the monomial basis
\end{lstlisting}}

Now we have a polynomial object to work with. A polynomial will always be seen as a formal sum of vectors: it cannot be factorized. If two polynomials are multiplied, the result will always be expanded as a sum.

{\tt \begin{lstlisting}
sage: pol * pol
x[2, 2, 4, 0] + 2*x[3, 4, 2, 0] + x[4, 6, 0, 0]
\end{lstlisting}}

\subsection{Apply operators}

We can now apply the \textit{divided differences} operators defined in Section \ref{DIFFDIV}.

{\tt \begin{lstlisting}
sage: pol = m[1,1,2] + m[2,3]; pol
x[1, 1, 2] + x[2, 3, 0]
sage: pol.divided_difference(2)
-x[1, 1, 1] + x[2, 1, 1] + x[2, 2, 0] + x[2, 0, 2]
sage: pol.divided_difference_isobar(2)
x[2, 1, 2] + x[2, 2, 1] + x[2, 3, 0] + x[2, 0, 3]
\end{lstlisting}}

\noindent
By default, the operator type is $A$, but we can also apply $B$, $C$, and $D$ operators. 

{\tt \begin{lstlisting}
sage: pol.divided_difference(2,"B")
x[1, -1, 2] + x[1, 0, 2] + x[2, 0, 0] + x[2, -3, 0] + x[2, -2, 0] + x[2, -1, 0] + x[2, 1, 0] + x[2, 2, 0]
sage: pol.divided_difference(2,"C")
x[1, 0, 2] + x[2, 0, 0] + x[2, -2, 0] + x[2, 2, 0]
sage: pol.divided_difference(2,"D")
x[0, 0, 0] + x[-2, -2, 0] + x[-1, -1, 0] + x[1, 1, 0] + x[1, 0, 2] + x[2, 2, 0] + x[0, -1, 2]
\end{lstlisting}}

We have seen in Section \ref{DIFFDIV} that type $B$, $C$, and $D$ operators were defined from the root systems of types $B_n$, $C_n$, and $D_n$, by the addition of the $n^{th}$ simple root and coroot. For each type, only one new operator was created which was the $n^{th}$ divided difference. But on the above example, applying the second divided difference on a polynomial in three variables gives a different result for types $B$, $C$, and $D$ than for type $A$. Even though groups have been used to give a definition of our operators, we have extended it to obtain a definition only depending on the polynomial. From the root systems of type $B_n$, $C_n$, and $D_n$, we had
\begin{align}
\partial_n^B &= \frac{1 - s_n^B}{x_n^{\frac{1}{2}} - x_{n}^{-\frac{1}{2}}} \\
\partial_n^C &= \frac{1 - s_n^C}{x_n - x_n^{-1}} \\
\partial_n^D &= \frac{1 - s_n^D}{x_{n-1}^{-1} - x_n}
\end{align}
So we just set
\begin{align}
\partial_i^B &= \frac{1 - s_i^B}{x_i^{\frac{1}{2}} - x_{i}^{-\frac{1}{2}}} \\
\partial_i^C &= \frac{1 - s_i^C}{x_i - x_i^{-1}} \\
\partial_i^D &= \frac{1 - s_i^D}{x_{i-1}^{-1} - x_i}
\end{align}
with $1 \leq i \leq n$ for types $B$ and $C$, and $ 2 \leq i \leq n$ for type $D$. These definitions allow us to study these operators for themselves and not relatively to their groups. As an example, this allows us to study the relations between $\partial_i^B$ and $\partial_{i+1}^B$. It does not make sense from a group point of view but it remains an interesting question. 

Of course, one may want to use the group based operators only and not the generalized ones. This is possible by using a different basis. The basis we have been using so far is called the \textit{Monomial basis} and is not related to any group. One can use a basis which is directly related to a root system called the \textit{Ambient space basis}. 

{\tt \begin{lstlisting}
sage: ma = A.ambient_space_basis("A"); ma
The ring of multivariate polynomials on x over Rational Field on the Ambient space basis of type A
sage: pol = ma[1,1,2] + ma[2,3]
sage: pol
x(1, 1, 2) + x(2, 3, 0)
sage: pol.parent()
The ring of multivariate polynomials on x over Rational Field with 3 variables on the Ambient space basis of type A
sage: pol.divided_difference(2)
-x(1, 1, 1) + x(2, 1, 1) + x(2, 2, 0) + x(2, 0, 2)
\end{lstlisting}}

\noindent
Note that \textit{Ambient Space basis} is very close to the \textit{Monomial basis} but the polynomial contains its type within its parent. It is related to a root system and only operators defined by this root system can be applied. 

{\tt \begin{lstlisting}
sage: mb = A.ambient_space_basis("B"); mb
The ring of multivariate polynomials on x over Rational Field on the Ambient space basis of type B
sage: pol = mb[1,1,2] + mb[2,3]
sage: pol.divided_difference(2)
-x(1, 1, 1) + x(2, 1, 1) + x(2, 2, 0) + x(2, 0, 2)
sage: pol.divided_difference(3)
x(1, 1, 0) + x(1, 1, -2) + x(1, 1, -1) + x(1, 1, 1)
\end{lstlisting}}

Conversions between the \textit{Monomial basis} and the \textit{Ambient space bases} can be carried out easily:

{\tt \begin{lstlisting}
sage: pol = m[1,1,2] + m[2,3]; pol
x[1, 1, 2] + x[2, 3, 0]
sage: pol.parent()
The ring of multivariate polynomials on x over Rational Field with 3 variables on the monomial basis
sage: pol = ma(pol); pol
x(1, 1, 2) + x(2, 3, 0)
sage: pol.parent()
The ring of multivariate polynomials on x over Rational Field with 3 variables on the Ambient space basis of type A
sage: pol = mb(pol); pol
x(1, 1, 2) + x(2, 3, 0)
sage: pol.parent()
The ring of multivariate polynomials on x over Rational Field with 3 variables on the Ambient space basis of type B
\end{lstlisting}}

\noindent
Even though the objects seem similar, one must always be careful with which basis one is working as it will impact the result of the operations as soon as operators are used. 

\subsection{Working with multi-bases}
\label{MULTIBASES}

We have already seen that our polynomials could be expressed on a different basis depending on which operations we wanted to make. But the \textit{Monomial basis} as well as the \textit{Ambient space bases} are just different versions of polynomials seen as sums of monomials. It is also possible to work with the linear bases we have defined in Section \ref{BASES}. Here is an example of the Schubert basis:

{\tt \begin{lstlisting}
sage: A = AbstractPolynomialRing(QQ)
sage: Schub = A.schubert_basis_on_vectors()
sage: Schub
The ring of multivariate polynomials on x over Rational Field on the Schubert basis of type A (indexed by vectors)
\end{lstlisting}}

\noindent
It can be used to create a Schubert polynomial and convert it to the monomial basis.

{\tt \begin{lstlisting}
sage: pol = Schub[1,2,2] + Schub[3,4]; pol
Y(1, 2, 2) + Y(3, 4, 0)
sage: pol.expand()
x(1, 2, 2) + x(2, 1, 2) + x(2, 2, 1) + x(3, 4, 0) + x(4, 3, 0)
sage: m(pol)
x[1, 2, 2] + x[2, 1, 2] + x[2, 2, 1] + x[3, 4, 0] + x[4, 3, 0]
sage: Schub(m[1,2,4] + m[2,3])
Y(1, 2, 4) - Y(1, 3, 3) - Y(1, 4, 2) - Y(2, 1, 4) + Y(2, 3, 0) + Y(2, 3, 2) + Y(2, 4, 1) + Y(3, 1, 3) - Y(3, 2, 0) - Y(3, 2, 2) - Y(4, 2, 1) + Y(5, 1, 1)
\end{lstlisting}}

\noindent
One can multiply Schubert polynomials together and the result will be given in the same basis. However the program is converting the two polynomials into the monomial basis to multiply them and then convert the result back into the Schubert basis.

{\tt \begin{lstlisting}
sage: pol1 = Schub[1,2,2] + Schub[3,4]
sage: pol2 = Schub[3,1,2]
sage: pol1 * pol2
Y(4, 3, 4) + Y(5, 2, 4) + Y(6, 5, 2) + Y(6, 6, 1) + Y(7, 4, 2) + Y(7, 5, 1)
\end{lstlisting}}

We have other bases implemented. Below is an example of the Key polynomials. One can convert directly from Schubert to Key polynomials without using the monomial basis manually (it is automatically done by the program):

{\tt \begin{lstlisting}
sage: K = A.demazure_basis_on_vectors();K
The ring of multivariate polynomials on x over Rational Field on the Demazure basis of type A (indexed by vectors)
sage: pol = K[2,1,4] + K[3,5,1];pol
K(2, 1, 4) + K(3, 5, 1)
sage: pol.expand()
x(2, 1, 4) + x(2, 2, 3) + x(2, 3, 2) + x(2, 4, 1) + x(3, 1, 3) + x(3, 2, 2) + x(3, 3, 1) + x(3, 5, 1) + x(4, 1, 2) + x(4, 2, 1) + x(4, 4, 1) + x(5, 3, 1)
sage: Schub(pol)
Y(2, 1, 4) + Y(3, 5, 1) - Y(5, 1, 1)
sage: K(m[1,2,4] + m[2,3])
K(1, 2, 4) - K(1, 3, 3) - K(1, 4, 2) - K(2, 1, 4) + K(2, 3, 0) + K(2, 3, 2) + K(2, 4, 1) + K(3, 1, 3) - K(3, 2, 0) - K(3, 2, 2) + K(4, 1, 2) - K(4, 2, 1)
sage: Khat = A.demazure_hat_basis_on_vectors()
sage: pol = Khat[2,1,4] + Khat[3,5,1];pol
^K(2, 1, 4) + ^K(3, 5, 1)
sage: pol.expand()
x(2, 1, 4) + x(2, 2, 3) + x(2, 3, 2) + x(3, 1, 3) + x(3, 2, 2) + x(3, 5, 1) + x(4, 4, 1)
sage: Schub(pol)
Y(2, 1, 4) - Y(2, 4, 1) + Y(3, 5, 1) - Y(4, 1, 2) + Y(4, 2, 1) - Y(5, 1, 1) - Y(5, 3, 1)
sage: Khat(m[1,2,4] + m[2,3])
^K(1, 2, 4) - ^K(1, 3, 3) + ^K(2, 3, 0) + ^K(2, 3, 2)
\end{lstlisting}}

\noindent
The key polynomials are also defined in type $B$, $C$, and $D$.

{\tt \begin{lstlisting}
sage: K = A.demazure_basis_on_vectors("B");K
The ring of multivariate polynomials on x over Rational Field on the Demazure basis of type B (indexed by vectors)
sage: pol = K[1,2,-2]
sage: pol
K(1, 2, -2)
sage: pol.expand()
x(1, 2, 0) + x(1, 2, -2) + x(1, 2, -1) + x(1, 2, 1) + x(1, 2, 2) + x(2, 1, 0) + x(2, 1, -2) + x(2, 1, -1) + x(2, 1, 1) + x(2, 1, 2) + x(2, 2, 0) + x(2, 2, -1) + x(2, 2, 1)
sage: pol = m[-2,1,1] + m[1,-1,1]; pol
x[-2, 1, 1] + x[1, -1, 1]
sage: K(pol)
K(0, 0, 0) + K(-2, 1, 1) - K(-1, 1, 1) - K(-1, 1, 2) - K(-1, 0, 1) - 2*K(1, 0, 0) - K(1, -2, 1) + K(1, -1, 0) + 2*K(1, -1, 1) + K(1, -1, 2) + K(1, 1, 0) - K(1, 1, -1) - 2*K(1, 0, 1) + K(0, 1, 1) + K(0, 0, 1)
\end{lstlisting}}

Back in type $A$, we also have the simple Grothendieck basis. We have two versions of it related by a change of variable explained in Section \ref{BASES} to avoid using negative exponents.

{\tt \begin{lstlisting}
sage: Grothn = A.grothendieck_negative_basis_on_vectors(); Grothn
The ring of multivariate polynomials on x over Rational Field on the Grothendieck basis of type A with negative exponents (indexed by vectors)
sage: pol = Grothn[1,2] + Grothn[2,2]; pol
G(1, 2) + G(2, 2)
sage: Grothp = A.grothendieck_positive_basis_on_vectors(); Grothp
The ring of multivariate polynomials on x over Rational Field on the Grothendieck basis of type A, with positive exponents (indexed by vectors) 
sage: pol.expand()
2*x(0, 0) + x(-2, 0) - x(-2, -1) - 3*x(-1, 0) - x(-1, -2) + 4*x(-1, -1) + x(0, -2) - 3*x(0, -1)
sage: pol = Grothp[1,2] + Grothp[2,2]; pol
G(1, 2) + G(2, 2)
sage: pol.expand()
x(1, 2) + x(2, 1)
sage: pol.expand().subs_var([(i,1-A.var(i)^(-1)) for i in xrange(1,3)])
2*x(0, 0) + x(-2, 0) - x(-2, -1) - 3*x(-1, 0) - x(-1, -2) + 4*x(-1, -1) + x(0, -2) - 3*x(0, -1)
\end{lstlisting}}

The last basis we have implemented are the non-symmetric Macdonald polynomials. In order to use it, one has to define a polynomial ring on a bigger field than $\mathbb{Q}$ to use variables in the coefficients.

{\tt \begin{lstlisting}
sage: var('t1 t2 q')
(t1, t2, q)
sage: K.<t1,t2,q> = QQ[]
sage: K = K.fraction_field()
sage: A = AbstractPolynomialRing(K);A
The abstract ring of multivariate polynomials on x over Fraction Field of Multivariate Polynomial Ring in t1, t2, q over Rational Field
sage: Mac = A.macdonald_basis_on_vectors()
sage: pol = Mac[1,2]; pol
M(1, 2)
sage: pol.expand()
t2^3*x(0, 0) + t2^2/q*x(1, 0) + ((t2*q+t2)/q^2)*x(1, 1) + 1/q^2*x(1, 2) + ((t2^2*q+t2^2)/q)*x(0, 1) + t2/q*x(0, 2)
sage: m = A.monomial_basis()
sage: Mac( m[1,1] )
(-t1*t2)*M(0, 0) + M(1, 0) + q*M(1, 1) + ((t1*t2*q^2+t2^2*q-t1*t2-t2^2)/(-t1*q-t2))*M(0, 1)
\end{lstlisting}}

\subsection{Define a new basis}

All our bases are defined with divided differences acting recursively on sums of monomials. If one needs to work with a new basis where objects are indexed by vectors, the only thing to be implemented is the rule converting one vector to a sum of monomials. The inverse conversion is automatically done if the basis is triangular. One can define one's own conversion function and so create a new basis. Below is an example to recreate the Schubert polynomials.

{\tt \begin{lstlisting}
sage: def schubert_on_basis(v, basis, call_back):
...     for i in xrange(len(v)-1):
...         if(v[i]<v[i+1]):
...             v[i], v[i+1] = v[i+1] + 1, v[i]
...             return call_back(v).divided_difference(i+1)
...     return basis(v)
\end{lstlisting}}

\noindent
Above is a definition of a recursive function that transforms a Schubert element into a sum of monomials. The principle is easy: we take the Schubert vector as an argument ($v$) and test if we find an index $i$ such that $v_i < v_{i+1}$. If we do, we compute the Schubert polynomial where $v_i$ and $v_{i+1}$ are switched and apply a divided difference operator. If $v$ is antidominant, then the result is $x^v$ that we get with $basis(v)$. The three parameters of this function are the vector $v$ corresponding to our Schubert element, a $basis$ parameter which is the monomial basis we convert to, and $call\_back$ that we use to call back our function (this ensures getting a cached method, 
i.e., things are not calculated twice). To create a new basis, one needs to write a function that takes these three arguments and returns the polynomial associated with the vector $v$. It can directly be written into the sage command line as above or in the notebook, or in a file attached to your session. The function will not directly be called by the user, it will be passed to a method and the program will send the right values to $basis$ and $call\_back$. If more arguments are needed, one just adds them this way:

{\tt \begin{lstlisting}
sage: def qt_schubert_on_basis(v, basis, call_back, q=1, t=1):
...     for i in xrange(len(v)-1):
...         if(v[i]<v[i+1]):
...             v[i], v[i+1] = v[i+1] + 1, v[i]
...             return q*1/t*call_back(v).divided_difference(i+1)
...     return basis(v)
\end{lstlisting}}

\noindent
Now that we have the function, we will pass it to our algebra to create a new basis:

{\tt \begin{lstlisting}
sage: A = AbstractPolynomialRing(QQ)
sage: myBasis = A.linear_basis_on_vectors("A","MySchub","Y",schubert_on_basis)
sage: pol = myBasis[2,1,3];pol
Y(2, 1, 3)
sage: pol.expand()
x(2, 1, 3) + x(2, 2, 2) + x(2, 3, 1) + x(3, 1, 2) + x(3, 2, 1) + x(4, 1, 1)
sage: myBasis(A.an_element())
Y(1, 2, 3) - Y(1, 3, 2) - Y(2, 1, 3) + Y(2, 3, 1) + Y(3, 1, 2) - Y(3, 2, 1) + Y(4, 1, 1)
\end{lstlisting}}

\noindent
This is a copy of the Schubert basis, and it works the same way as the previous bases we have seen in Section \ref{MULTIBASES}. Below is an example with a parametrized function:

{\tt \begin{lstlisting}
sage: var('q t')                                                                                         
(q, t)
sage: K.<q,t> = QQ[]
sage: K = K.fraction_field()
sage: A = AbstractPolynomialRing(K)
sage: qtSchubertBasis = A.linear_basis_on_vectors("A","qtSchub","YQ",qt_schubert_on_basis,(("q",q),("t",t)) )
sage: pol = qtSchubertBasis[1,2,3]; pol
YQ(1, 2, 3)
sage: pol.expand()
q^3/t^3*x(1, 2, 3) + q^3/t^3*x(1, 3, 2) + q^3/t^3*x(2, 1, 3) + 2*q^3/t^3*x(2, 2, 2) + q^3/t^3*x(2, 3, 1) + q^3/t^3*x(3, 1, 2) + q^3/t^3*x(3, 2, 1)
\end{lstlisting}}

\noindent
The extra parameters are sent by a tuple of couples \textit{(parameter name, parameter value)} to the main algebra that will create the new basis. 

\subsection{Double set of variables}

Our program also contains another algebra to work with a double set of variables.

{\tt \begin{lstlisting}
sage: D = DoubleAbstractPolynomialRing(QQ); D
The abstract ring of multivariate polynomials on x over The abstract ring of multivariate polynomials on y over Rational Field
\end{lstlisting}}

\noindent
One can see that the double algebra is the algebra of multivariate polynomials in the $x$ variables with the multivariate polynomials in the $y$ variables as coefficients.

{\tt \begin{lstlisting}
sage: D.an_element()
y[0]*x[0, 0, 0] + 2*y[0]*x[1, 0, 0] + y[0]*x[1, 2, 3] + 3*y[0]*x[2, 0, 0]
\end{lstlisting}}

\noindent
One can specify which bases to use in the $x$ variables and which bases to use in the $y$ variables.

{\tt \begin{lstlisting}
sage: Dx = D 
sage: Dy = D.base_ring()
sage: Schubx = Dx.schubert_basis_on_vectors()
sage: Schuby = Dy.schubert_basis_on_vectors()
sage: pol = Schuby[2,1,3] * Schubx[1,1,2]
sage: pol
(Yy(2,1,3))*Yx(1, 1, 2)
\end{lstlisting}}

\noindent
The \textit{expand} function or all direct conversions are done in the $x$ variables.

{\tt \begin{lstlisting}
sage: pol.expand()
(Yy(2,1,3))*x(1, 1, 2) + (Yy(2,1,3))*x(1, 2, 1) + (Yy(2,1,3))*x(2, 1, 1)
sage: mx = Dx.monomial_basis()
sage: mx(pol)
(Yy(2,1,3))*x[1, 1, 2] + (Yy(2,1,3))*x[1, 2, 1] + (Yy(2,1,3))*x[2, 1, 1]
\end{lstlisting}}

\noindent
But, of course, one can also easily change the basis for the $y$ variables.

{\tt \begin{lstlisting}
sage: my = Dy.monomial_basis()
sage: pol.change_coeffs_bases(my)
(y[2,1,3]+y[2,2,2]+y[2,3,1]+y[3,1,2]+y[3,2,1]+y[4,1,1])*Yx(1, 1, 2)
sage: pol = mx(pol); pol
(Yy(2,1,3))*x[1, 1, 2] + (Yy(2,1,3))*x[1, 2, 1] + (Yy(2,1,3))*x[2, 1, 1]
sage: pol.change_coeffs_bases(my)
(y[2,1,3]+y[2,2,2]+y[2,3,1]+y[3,1,2]+y[3,2,1]+y[4,1,1])*x[1, 1, 2] + (y[2,1,3]+y[2,2,2]+y[2,3,1]+y[3,1,2]+y[3,2,1]+y[4,1,1])*x[1, 2, 1] + (y[2,1,3]+y[2,2,2]+y[2,3,1]+y[3,1,2]+y[3,2,1]+y[4,1,1])*x[2, 1, 1]
\end{lstlisting}}

\noindent
One can also change the role of variables between the main ones and the coefficients.

{\tt \begin{lstlisting}
sage: pol
(Yy(2,1,3))*x[1, 1, 2] + (Yy(2,1,3))*x[1, 2, 1] + (Yy(2,1,3))*x[2, 1, 1]
sage: pol.swap_coeffs_elements()
(x[1,1,2]+x[1,2,1]+x[2,1,1])*Yy(2, 1, 3)
\end{lstlisting}}

So we have seen that we can use our previous bases on a double set of variables. But we also have specific bases that only work with a double set of variables. Let us see the double Schubert polynomials and double Grothendieck polynomials, the way they were defined in Section \ref{BASES}.

{\tt \begin{lstlisting}
sage: DoubleSchub = D.double_schubert_basis_on_vectors(); DoubleSchub
The ring of multivariate polynomials on x over The abstract ring of multivariate polynomials on y over Rational Field on the Double Schubert basis of type A (indexed by vectors)
sage: pol = DoubleSchub[1,2]; pol
y[0]*YY(1, 2)
sage: pol.expand()
(-y(2,1,0)-y(2,0,1))*x(0, 0) + (y(1,1,0)+y(1,0,1)+y(2,0,0))*x(1, 0) + (-2*y(1,0,0)-y(0,1,0)-y(0,0,1))*x(1, 1) + y[0]*x(1, 2) + (-y[1])*x(2, 0) + y[0]*x(2, 1) + (y(1,1,0)+y(1,0,1)+y(2,0,0))*x(0, 1) + (-y[1])*x(0, 2)
sage: DGroth = D.double_grothendieck_basis_on_vectors(); DGroth
The ring of multivariate polynomials on x over The abstract ring of multivariate polynomials on y over Rational Field on the Double Grothendieck basis of type A (indexed by vectors)
sage: pol =DGroth[1,2]; pol
y[0]*GG(1, 2)
sage: pol.expand()
y[0]*x(0, 0) + (-y(2,1,1))*x(-2, -2) + (y(1,1,1))*x(-2, -1) + (-y[1])*x(-1, 0) + (y(1,1,1))*x(-1, -2) + (y(2,0,0)-y(0,1,1))*x(-1, -1) + (-y[1])*x(0, -1)
\end{lstlisting}}

\section{Some advanced applications}

\subsection{Projective degrees of Schubert varieties}

In his Ph.D. thesis, Veigneau \cite{VEIGN} presents an application of the software ACE by computing the projective degree of Schubert varieties. We can implement this same function with our patch on Sage-Combinat. The projective degree $d(X)$ of a sub-variety $X \subset \mathbb{P}^M$ of codimension $k$ is the number of intersections between $X$ and a generic hyperplane of dimension $k$. For $\sigma$ a permutation of size $n$ and $X_\sigma$ a Schubert sub-variety of the flag variety $\mathcal{F}(\mathbb{C}^n)$ embedded in $\mathbb{P}^M$ by the Pl\"ucker embedding (with $M = 2^N-1$ where $N=\frac{n(n-1)}{2}$ is the dimension of $\mathcal{F}(\mathbb{C}^n)$), $d(X_\sigma)$ is a coefficient in a product in the Schubert basis. More precisely, the first Chern class of the tautologic invertible vectorial fiber of $\mathbb{P}^M$ is $h = (n-1)x_1 + (n-2)x_2 + \ldots + x_{n-1}$ and $d(X_\sigma)$ is the coefficient of $Y_{n-1,n-2,\ldots,0}$ in  $h^{N-\ell(\sigma)} Y_v$, where $Y_v$ is the Schubert polynomial indexed by $v$, the Lehmer code of $\sigma$ \cite{CHERN}. The following function computes these degrees:

{\tt \begin{lstlisting}
def proj_deg(perm):
    n = len(perm)
    d = n*(n-1)/2 - perm.length()
    
    # we create the polynomial ring and the bases
    A = AbstractPolynomialRing(QQ)
    Schub = A.schubert_basis_on_vectors()
    
    # we compute the product
    h = sum( [(n-i) * A.var(i) for i in xrange(1,n)])
    res = Schub( h**d * Schub(perm.to_lehmer_code()))
    
    # we look for the right coefficient
    for (key, coeff) in res:
        if ([key[i] for i in xrange(n)] == [n-i for i in xrange(1,n+1)]):
            return coeff
            
    return 0
\end{lstlisting}} 

\noindent
One can also compute the product and directly read the result:

{\tt \begin{lstlisting}
sage: A = AbstractPolynomialRing(QQ)
sage: m = A.monomial_basis()
sage: Schub =  A.schubert_basis_on_vectors()
sage: Schub( (3*m[1] + 2*m[0,1] + m[0,0,1])^4 * Schub[1,0,1,0])
8*Y(1, 1, 4, 0) + 23*Y(1, 2, 3, 0) + 24*Y(1, 3, 2, 0) + 39*Y(1, 4, 1, 0) + 15*Y(1, 5, 0, 0) + Y(1, 0, 5, 0) + 48*Y(2, 1, 3, 0) + 101*Y(2, 2, 2, 0) + 117*Y(2, 3, 1, 0) + 84*Y(2, 4, 0, 0) + 12*Y(2, 0, 4, 0) + 173*Y(3, 1, 2, 0) + 78*Y(3, 2, 1, 0) + 147*Y(3, 3, 0, 0) + 53*Y(3, 0, 3, 0) + 283*Y(4, 1, 1, 0) + 171*Y(4, 2, 0, 0) + 96*Y(4, 0, 2, 0) + 93*Y(5, 1, 0, 0) + 176*Y(5, 0, 1, 0) + 80*Y(6, 0, 0, 0)
sage: proj_deg(Permutation([2,1,4,3]))
78
\end{lstlisting}}

\noindent
We can use our function to compute the degree for all permutations of size 4:

{\tt \begin{lstlisting}
sage: degrees = {}
sage: for perm in Permutations(4):
....:     degrees[perm] = proj_deg(perm)
....:     
sage: degrees
{[2, 1, 4, 3]: 78, [1, 3, 4, 2]: 48, [3, 2, 4, 1]: 3, [3, 1, 2, 4]: 48, [4, 2, 1, 3]: 3, [1, 4, 2, 3]: 46, [3, 2, 1, 4]: 16, [4, 1, 3, 2]: 3, [2, 3, 4, 1]: 6, [3, 4, 2, 1]: 1, [1, 2, 3, 4]: 720, [1, 3, 2, 4]: 280, [2, 4, 3, 1]: 3, [2, 3, 1, 4]: 46, [3, 4, 1, 2]: 2, [4, 2, 3, 1]: 1, [1, 4, 3, 2]: 16, [4, 1, 2, 3]: 6, [2, 4, 1, 3]: 12, [4, 3, 1, 2]: 1, [4, 3, 2, 1]: 1, [3, 1, 4, 2]: 14, [2, 1, 3, 4]: 220, [1, 2, 4, 3]: 220}
\end{lstlisting}}

\subsection{Determinants of Schur functions}

Grassmannian Schubert polynomials are the Schubert polynomials indexed by vectors $v$ such that $v_1 \leq v_2 \leq \ldots \leq v_n$. They are symmetric functions in $x_1, \ldots, x_n$. In a single set of variables (i.e., specializing $y$ to 0), Grassmannian Schubert polynomials are equal to Schur functions. More precisely, the transition matrix between double Grassmannian Schubert polynomials and Schur functions is unitriangular. 

{\tt \begin{lstlisting}
sage: A = AbstractPolynomialRing(QQ)
sage: Schub =  A.schubert_basis_on_vectors()
sage: pol = Schub[1,2]
sage: pol.expand()
x(1, 2) + x(2, 1)
sage: 
sage: D = DoubleAbstractPolynomialRing(QQ)
sage: DSChub = D.double_schubert_basis_on_vectors()
sage: pol = DSchub[1,2]
sage: pol
y[0]*YY(1, 2)
sage: Schub = D.schubert_basis_on_vectors()
sage: Schub(pol)
y[0]*Yx(1, 2) + (-y(2,1,0)-y(2,0,1))*Yx(0, 0) + (-y(1,0,0)-y(0,1,0)-y(0,0,1))*Yx(1, 1) + (y(1,1,0)+y(1,0,1)+y(2,0,0))*Yx(0, 1) + (-y[1])*Yx(0, 2)
\end{lstlisting}}

\noindent
This allows us to compute determinants of Schur functions by replacing these by Schubert polynomials and specializing arbitrarily the $y$ variables. For example, we can compute
\begin{equation}
\vert s_\mu(A) \vert_{\mu \subseteq 11},
\end{equation}
where $A \in [ \lbrace x_1, x_2 \rbrace, \lbrace x_1, x_3 \rbrace, \lbrace x_2, x_3 \rbrace ]$, and prove that it is equal to $\prod_{j>i}(x_j-x_i)$. First, we replace $s_\mu$ by 
\begin{equation}
\vert Y_u(A,y) \vert_{u = 00, 01, 11}
\end{equation}
and specialize $y_1=x_1$, $y_2=x_2$, in which case the determinant becomes
\begin{equation}
\left|
\begin{array}{lll}
1 & 1 & 1 \\
0 & x_3 - x_2 & x_3 - x_1 \\
0 & 0 & (x_3-x_1)(x_2-x_1)
\end{array}
\right|
\end{equation}
and gives the result. The following function computes the above Matrix:

{\tt \begin{lstlisting}

def compute_matrix(variables, alphabet, indices):
    
    n = len(indices)
    
    #Initial definitions
    K = PolynomialRing(QQ,[var(v) for v in variables])
    K = K.fraction_field()
    D = DoubleAbstractPolynomialRing(K)
    DSchub = D.double_schubert_basis_on_vectors()
    
    result_matrix = []
    
    for u in indices:
        line = []
        
        #the expansion on the double schubert will allow us to compute the result
        pu = DSchub(u).expand()
        
        for a in alphabet:
        	#we apply our polynomial on alphabets and specialize the y (this should be improved on further versions)
            pol = pu.subs_var( [(i,K(a[i])) for i in xrange(len(a))])
            pol = pol.swap_coeffs_elements()
            pol = pol.subs_var( [(i,K(variables[i])) for i in xrange(pol.nb_variables()) ] )
            if(pol ==0):
                coeff = 0
            else:
                coeff = list(list(pol)[0][1])[0][1]
            line.append(coeff)
        result_matrix.append(line)
    
    return Matrix(K,result_matrix)
\end{lstlisting}}

\noindent
In Sage:

{\tt \begin{lstlisting}
sage: variables = ("x1", "x2", "x3")
sage: alphabet = [["x1","x2"],["x1","x3"],["x2","x3"]]
sage: indices = [[0,0],[0,1],[1,1]]
sage: 
sage: res = compute_matrix(variables, alphabet, indices)
sage: res
[  1         1               1                           ]
[  0         -x2 + x3        -x1 + x3                    ]
[  0         0               x1^2 - x1*x2 - x1*x3 + x2*x3]
sage: det = res.determinant()
sage: det
-x1^2*x2 + x1*x2^2 + x1^2*x3 - x2^2*x3 - x1*x3^2 + x2*x3^2
sage: factor(det)
(x2 - x3) * (-x1 + x2) * (x1 - x3)
\end{lstlisting}}


\end{document}